\documentclass[12pt]
{amsart}

\usepackage{amsfonts,amsmath,amsthm,amssymb}
\usepackage[mathscr]{euscript}
\usepackage{graphicx}
\usepackage{multicol}

\newcommand{\C}{{\mathbb C}}

\newcommand{\disk}{{\mathbb D}}

\newcommand{\ovl}[1]{\overline{#1}}

\newcommand{\eps}{\varepsilon}
\newcommand{\Cbar}{\ovl{\C}}

\renewcommand{\phi}{\varphi}

\newcommand{\reminder}[1]{\textsf{#1}}
\newcommand{\sergey}[1]{\textsl{#1}}	
\newcommand{\hide}[1]{}

\newtheorem{lemma}{Lemma}
\newtheorem{theorem}[lemma]{Theorem}

\newtheorem{conjecture}{Conjecture}

\newcommand{\pack}[1]{\texttt{#1}}

\theoremstyle{definition}

\allowdisplaybreaks

\title[Root finding by dynamical systems]{
Finding polynomial roots by dynamical systems --- a case study}

\author[S. Shemyakov]{Sergey Shemyakov}
\address{Institut de Mathématiques (UMR CNRS7373) Campus de Luminy 163 avenue de Luminy - Case 907 13288 Marseille 9}
\email{sergeyshemyakov@gmx.de}
\author[et al.]{Roman Chernov}
\email{roman1chernov@gmail.com}
\author[]{Dzmitry Rumiantsau}
\email{d.rumiantsau@jacobs-university.de}
\author[]{Dierk Schleicher}
\email{dierk.schleicher@gmx.de}
\author[]{Simon Schmitt}
\author[]{Anton Shemyakov}
\email{a.shemyakov@jacobs-university.de}

\begin{document}

\begin{abstract}
We investigate two well known dynamical systems that are designed to find roots of univariate polynomials by iteration: the methods known by Newton and by Ehrlich--Aberth. Both are known to have found all roots of high degree polynomials with good complexity. Our goal is to determine in which cases which of the two algorithms is more efficient. We come to the conclusion that Newton is faster when the polynomials are given by recursion so they can be evaluated in logarithmic time with respect to the degree, or when all the roots are all near the boundary of their convex hull. Conversely, Ehrlich--Aberth has the advantage when no fast evaluation of the polynomials is available, and when roots are in the interior of the convex hull of other roots. 
\end{abstract}

\maketitle

\section{Introduction}

Among the most fundamental and classical results in mathematics is the fact that every univariate polynomial of degree $d$ over $\C$ splits into $d$ linear factors (i.e., $\C$ is algebraically closed), and that there is no closed formula to determine their roots in terms of radicals. Therefore, these roots have to be approximated numerically, usually in terms of an iterated process that can be viewed as a very natural dynamical system.

We investigate two such dynamical systems that are particularly prominent as root finders: Newton's method and the Ehrlich--Aberth-method. Both are well known as efficient root finders that have a good track record for finding all roots of rather large degrees, and the question arises which of these is more efficient. Of course, there are numerous other root finding procedures, including eigenvalue methods or Weierstrass' iteration method. Our choice is motivated by the fact that the Ehrlich--Aberth-method is underlying one of the best established practical root finding software packages, \pack{MPSolve}, while in our own work Newton's method has turned out as particularly successful. Moreover, these two root finders are conceptually related in an interesting way: Newton approximates one root at a time, while Ehrlich--Aberth may be seen as a system of $d$ orbits that try to approximate all the $d$ roots at the same time, and these $d$ orbits can be viewed as $d$ Newton orbits that are synchronized in the sense that each of them is aware of where the others are: they are ``attracted'' by the roots in the same sense as for Newton's method, but repelled from each other so that they avoid approximating the same (simple) root.

We feel it would be quite interesting to extend this comparison also to other root finding procedures. 

It is clear that no continuous deterministic root finding system can converge to roots for every set of initial conditions: the domains of convergence to some set of roots is open, and the complement must contain at least the boundaries of the basins. The second best desired property of a root finder is that it be \emph{generally convergent}: this means it should converge to a root for an open and dense subset of all possible starting conditions (starting points for Newton, starting vectors for Ehrlich--Aberth). This property is difficult to establish: Newton's method is known not to be generally convergent, and for Ehrlich--Aberth it is an open question. 

The structure of this paper is as follows: in Section~\ref{Sec:NewtonEhrlichAberth} we describe the iteration steps of Newton and of Ehrlich--Aberth, trying to bring out their conceptual similarities. Relevant and interesting properties of both methods are then described, respectively, in Sections~\ref{Sec:Ehrlich-Aberth} and \ref{Sec:Newton}. Our key results are presented in Section~\ref{Sec:Experiments}: we describe several families of polynomials with different properties, and all with a wide range of degrees, and measure the performance of Newton's and Ehrlich--Aberth's methods in systematic experiments. The final Sections~\ref{Sec:Conclusion} contains a discussion of the results and the complexity of the method, as well possible  conclusions and ideas for improvements that might inspire further research.

\subsection*{Acknowledgements}
This project has been inspired by several of our friends and colleagues. It owes a lot to many interesting discussions with Dario Bini and Leonardo Robol in Pisa, as well as with Victor Pan in New York, as well as John Hubbard in Cornell. It builds on many ideas and discussions with the members of our research group, in particular Robin Stoll, Marvin Randig, and Bernhard Reinke. In addition, we are most grateful to the two referees for their helpful comments and suggestions.

We gratefully acknowledge support through the Advanced Grant HOLOGRAM of the  European Research Council.

\section{The root finding methods by Newton and Ehrlich--Aberth in comparison}
\label{Sec:NewtonEhrlichAberth}

The setting in all cases is the same: we investigate a monic polynomial $p(z)=(z-\alpha_1)\dots(z-\alpha_d)$ of degree $d$, with $d$ roots $\alpha_1$,\dots,$\alpha_d$; nothing will be lost by assuming that $p$ is monic.

\subsection{Newton's method}

Newton's well known root finding method is not an algorithm but a heuristic: given a polynomial $p$ (or more generally a differentiable mapping) and a guess $z$ for a possible root, then a ``suggested improvement'' for the guess is
\begin{equation}
z'=N_p(z):=z-p(z)/p'(z)
\;.
\label{Eq:NewtonIteration}
\end{equation}
If $z$ is already a root, i.e.\ $p(z)=0$, then this point remains fixed: roots of $p$ are fixed points of $N_p$. In general, $N_p$ is a rational map of the same degree as the number of \emph{distinct} roots of $p$ (multiple roots of $p$ are accounted for in the degree of $p$ but not of $N_p$).

\begin{lemma}[Local convergence of Newton's method]
Every simple root $\alpha$ has a neighborhood on which convergence of Newton's method is quadratic: $|N_p(z)\alpha|<C |z-\alpha|^2$. For multiple roots the convergence is only linear: if $\alpha$ is a root with multiplicity $k\ge 2$, then $|N_p(z)-\alpha|/|z-\alpha|\longrightarrow (k-1)/k$ as $z\to\alpha$.
\end{lemma}

\looseness+1
There is an electrostatic interpretation of Newton's method: place protons at the locations of the roots of $p$; then the electrostatic force felt by an electron at a point $z\in\C$ is the sum of the attractive forces of all the roots, and each of these forces is in the direction to the root and with strength inversely proportional to the distance. Except for constants of physical units, the combined force equals $-\sum_i \ovl{1/(z-\alpha_i)}=-\ovl{p'(z)/p(z)}$ (where the overbar denotes complex conjugation). For Newton's method, the image $z'$ of $z$ is the unique point for which this sum is equal to $-\ovl{1/(z-z')}$: it places the point $z'$ so that a single charge at $z'$ has the same effect as the superposition of all $d$ charges at the roots. In other words, it moves $z$ in the direction of the net force of all the roots, but at a distance inversely proportional to the net force: the smaller this net force, the further away all the roots, and the further the point $z$ has to move in the direction of the origin of the force.

\hide{
\reminder{Mention here that Newton is not globally convergent?}
}

\goodbreak

\subsection{The Ehrlich--Aberth-method}

This method does not approximate one root at a time, but tries to find all $d$ roots at once; hence it is an iteration on $\C^d$ (or rather an open subset thereof).

Consider a vector $(z_1,\dots,z_d)$ of initial guesses. The Ehrlich--Aberth-iteration consists of the map $(z_1,\dots,z_k)\mapsto (z'_1,\dots,z'_k) $ given by
\begin{equation}
z'_k:=z_k-\frac{\frac{p(z_k)}{p'(z_k)}}{1-\frac{p(z_k)}{p'(z_k)}\cdot \sum_{i\neq k}\frac{1}{z_k-z_i}}
\;.
\label{Eq:EhrlichAberthIteration}
\end{equation}

Observe that if some $z_k$ is equal to a root of $p$, i.e.\ $p(z_k)=0$, then $z'_k=z_k$ so this coordinate remains fixed. There is a heuristic interpretation for this formula: each of the $d$ component variables $z_k$ ``thinks'' that all the other components are already equal to a root (``everyone else knows what they are doing, only I am wrong''), and tries to use the standard Newton update $z\mapsto z-f(z)/f'(z)$ for the function $f(z)=p(z)/\prod_{i\neq k}(z-z_k)$, which would be the linear function $(z-z_k)$ if all the other approximations were correct. To see this, observe that
\begin{align*}
f'(z)&=\frac{p'(z)\prod_{i\neq k}(z-z_k) - p(z) \prod_{i\neq k}(z-z_k)\sum_{i}\frac{1}{z_k-z_i}}{(\prod_{i\neq k}(z-z_k))^2}
\\
&= \frac{p'(z)-p(z)\sum_{i}\frac{1}{z_k-z_i}}{\prod_{i\neq k}(z-z_k)}
\end{align*}
and hence
\[
\frac{f(z_k)}{f'(z_k)}= \frac{p(z_k)}{p'(z_k)-p(z_k)\sum_{i}\frac{1}{z_k-z_i}}
= \frac{\frac{p(z_k)}{p'(z_k)}}{1-\frac{p(z_k)}{p'(z_k)}\cdot \sum_{i\neq k}\frac{1}{z_k-z_i}}
\]
as claimed above.

\begin{lemma}[Local convergence of Ehrlich--Aberth method]
\label{Lem:LocalConvergenceEhrlichAberth}
Every vector consisting of the $d$ roots of $p$ has a neighborhood in $\C^d$ on which the Ehrlich--Aberth method converges to this solution vector. This convergence is cubic when all roots of $p$ are simple, and linear otherwise.
\end{lemma}



This method has an electrostatic interpretation similar as Newton's method: once again, we place positive unit charges at the positions $\alpha_1,\dots,\alpha_d$ of the $d$ roots in $\C$. This time, there are $d$ approximations to the roots at the positions $z_1,\dots,z_k$, each coming with a negative unit charge. The $d$ approximations to the roots are as much attracted to the actual roots as they are repelled by each other, so they have a natural tendency to find different roots; in particular, if one approximation has already ``found'' a root, then the positive and negative charges cancel exactly.

\begin{lemma}[Electrostatic interpretation of Ehrlich--Aberth-method]
In the Ehrlich--Aberth-method, each root moves in the direction of the resulting net charge, and a distance that is inversely proportional to the strength of this charge.
\end{lemma}
\begin{proof}
The following computation starts with the Ehrlich--Aberth iteration formula end terminates with the electrostatic interpretation:
\begin{align*}
z'_k:=z_k-\frac{\frac{p(z_k)}{p'(z_k)}}{1-\frac{p(z_k)}{p'(z_k)}\cdot \sum_{i\neq k}\frac{1}{z_k-z_i}}
=z_k-\frac{1}{ \frac{p'(z_k)-p(z_k)\sum_{i\neq k}\frac{1}{z_k-z_i}}{p(z_k)}}\\
=z_k-\frac{1}{\frac{p'(z_k)}{p(z_k)}-\sum_{i\neq k}\frac{1}{z_k-z_i} }
= z_k-\frac{1}{\sum_i\frac{1}{z_k-\alpha_i}-\sum_{i\neq k}\frac{1}{z_k-z_i}  }
\;.
\end{align*}
So indeed, the formula is the same as Newton's method, except that all the other $d-1$ components in the vector of guesses are added as ``roots with negative charge''.
\end{proof}

If one runs Newton's method on $d$ approximations separately, one can give it a similar electrostatic interpretation in the sense that the $d$ approximations have infinitesimally small test charges so that they see the roots, but do not interact with each other. This can be seen as a conceptual advantage of Ehrlich--Aberth over Newton: the $d$ approximations are synchronized with each other. However, currently there does not seem to be theory known that exploits this fact towards an understanding of the global dynamics of the Ehrlich--Aberth-method. In particular, it is not known whether it is globally convergent; see Section~\ref{Sub:EhrlichAberthProperties}.

One key difference between the Newton and Ehrlich--Aberth-methods will be relevant below: in order to evaluate Newton's method, it suffices to have a way of evaluating $p'/p$. In certain cases, for instance  when $p$ is given by iteration or other efficient types of recursion, this may have much better complexity than $O(d)$: in the cases explored here, this complexity was $O(\log d)$. Linear recursions for $p$ do not offer improvements here, but sparse polynomials might (compare for instance \cite{BiniFiorentino}). 

 Even though one can use the same efficient evaluation of $p'/p$ in Ehrlich--Aberth, this does not improve the overall complexity of the method. This explains why specifically in such cases, the performance of Newton is much better.

We would like to mention that for both methods approximate improvements are indeed possible that might lead to substantially faster evaluations at least in appropriate cases: these are discussed in Section~\ref{Sub:Improvements}.


\goodbreak

\subsection{Machine precision and stopping criteria}  

Like for all numerical experiments, it is important to keep in mind the capabilities and limitations of the computing system. Our experiments were all performed with floating point numbers with double precision. In practice, we made the ``idealistic'' assumption that all our calculations were exact. They are not, of course, but we have reasons to believe that this is not a relevant issue in practice for our polynomials, even for large degrees: both numerical methods are known to be intrinsically stable with self-correcting errors, and one can often verify a posteriori  that indeed all roots have been found (see for instance in \cite[Section~2]{NewtonRobin1}, where an explicit worst-case estimate on error bounds for Newton's method is computed, and \cite{NewtonRobin2} for such an a-posteriori verification for degrees  up to $2^{30}$). In fact, in many cases (including all our polynomials) the mutual distance between roots is much larger compared to available precision of computation, so all roots have numerically distinguishable domains of quadratic convergence. In each case it remains an a-posteriori verification that the conceivable problem of insufficient machine precision does not lead to missed roots.

A final remark concerns the stopping criterion: a frequent condition that a root has been found is when the value of the polynomial is smaller than a given threshold. A different criterion is that a particular approximation $z_k$ is $\delta$-close to an actual root $\alpha_i$, with a precise bound on $\delta$. This can be upgraded to the requirement that all roots of $p$ have been found if there exists a vector $(z_1,\dots,z_d)\in\C^d$ with a guarantee that, up to permutation, $|z_i-\alpha_i|<\delta$ for all $i\in\{1\dots,d\}$. That is the stopping criterion that we used for Ehrlich--Aberth implementation. 

The question of whether some root has higher multiplicity is irrelevant (and numerically usually not a valid question): what matters is that several $z_i$ are close to each other, whether or not the approximated roots are multiple or only nearby. As we mentioned before, all the polynomials considered in our experiments have well-separated roots.

The actual precision of the computed roots is not of great relevance because of quadratic convergence, provided the roots are simple: once the roots have been found in the sense that they are separated then any desired precision $\eps$ can be achieved in time of about $O(d \log|\log\eps|)$ (within the limits of machine precision). For both methods, by far most of the computation time is spent on separating the roots. However if one wants to estimate the error from the perspective of numerical analysis, the paper \cite{Proinov} gives some useful a priori and a posteriori bounds for Newton's method.

\goodbreak

\section{The Ehrlich--Aberth-method}
\label{Sec:Ehrlich-Aberth}

\subsection{Properties of the method; general convergence}
\label{Sub:EhrlichAberthProperties}

The key feature of this method is the recursive step, given in \eqref{Eq:EhrlichAberthIteration}, that takes one vector of $d$ complex numbers, viewed as initial guesses, and computes a new vector of supposedly improved approximations of the $d$ roots of a given degree $d$ polynomial $p$. In order to upgrade this iteration step to an algorithm, one needs to specify the vector of initial guesses (starting points), as well as a stopping criterion. Moreover, at least in principle there is the possibility that the iteration fails to find the roots of the polynomials, and this needs to be detected. 

While the local convergence of the Ehrlich--Aberth-method in a neighborhood of the roots is understood (Lemma~\ref{Lem:LocalConvergenceEhrlichAberth}), we are not aware of any theory about global convergence properties. The method cannot converge in all cases. An obvious case where it fails is when $p$ is a real polynomial with non-real roots, and the vector of starting points is entirely real. Then by symmetry all subsequent iterates will be entirely real, so they cannot converge to the roots. There are similar symmetric situations that prevent convergence to roots. 

More conceptually, convergence must fail on larger sets of starting points: the set of initial vectors $(z_1,\dots,z_d)\in\C^d$ from which the iteration converges to a vector of roots in any particular order is open, and it is non-empty because it includes a neighborhood of the root vector. Since this is so for every ordered vector of the $d$ roots, it follows that the subset of $\C^d$ from which convergence occurs is a finite union of disjoint open sets, more than one, and this cannot be all of $\C^d$: therefore convergence must fail on all the boundaries, which must be large enough to separate open sets in $\C^d$, so it must fail on a set of topological dimension at least $2d-1$.

In practice, failure to converge seems rare, and apparently has not been observed except in the case of symmetries, despite extensive experience in implementations such as MPSolve; see \cite{BiniAberth,RobolMaster}. It may well be that, from a dynamical systems point of view, the locus of non-convergence is unstable: any small perturbation may lead to iteration away from this locus and convergence towards the roots; but for all we know this is an open question. There are periodic orbits: these are finitely many points in $\C^d$ that are permuted by the iteration. The way how  Newton's method fails to be generally convergent (see Section~\ref{Sec:Newton}) is because some of its periodic points may become attracting, so they attract a whole open set of starting points, and this happens for all non-trivial cases: even for degree $3$ there are  polynomials that have attracting cycles of any  period. By analogy, one may suspect that the Ehrlich--Aberth-method also has periodic cycles that are attracting. In fact, this  has been the intuition of several people in dynamical systems; however, in numerous experiments in particular using the \texttt{MPSolve} implementation no attracting cycles have ever been observed.

At this point a digression to the Weierstrass method may be interesting, even though it is not the focus of the current paper. Like the Ehrlich--Aberth-method, it is an iteration in $\C^d$ (undefined at certain points where two or more coordinates coincide). Both methods are known as reliable root finders that seem to find roots generically and efficiently, but both are lacking global theory. 

For the Weierstrass method, there are recent results in \cite{BernhardMichaelWeierstrass}. Most importantly, it has been discovered that \emph{the Weierstrass is not generally  convergent}: there is an open subset  of the space of polynomials of degree $3$ or higher that has an open set of  starting vectors that fails to converge to any  roots. More precisely, there are explicit cubic polynomials such as $z\mapsto z^3+z+180$ for which the Weierstrass method has periodic points of period $4$ that are attracting. This result seems to go against  the expectation of people working in numerical analysis: but it is what had been suspected in analogy to  Newton's method by people in dynamical systems, in particular by Smale.

The result in \cite{BernhardMichaelWeierstrass} comes out of systematic investigation of periodic points of small degrees and periods. It is shown that for cubic polynomials, period $4$ is minimal: all periodic points of period $2$ or $3$ are non-attracting because explicit equations are found for the sums of the eigenvalues at periodic points, and these equations are not compatible with the condition that all eigenvalues are in $\disk$ as required for attracting periodic points. Another unexpected property that was discovered is that there are starting points for which the Weierstrass iteration is defined at all times but for which the dynamics tends to $\infty$; this is so for (almost) every polynomial of degree $3$ or higher.

This situation is less clear for the Ehrlich--Aberth-method. This method seems similar to the Weierstrass method, but it is more complicated: the equations themselves are algebraically more complicated, and Weierstrass has the advantage that it has an invariant subspace of codimension $1$ where the interesting dynamics takes place. Therefore, the analysis in \cite{BernhardMichaelWeierstrass} (which even for Weierstrass very quickly developes high complexity that limits the cases that can be treated) can cover only basic cases for Ehrlich--Abert, and it is not clear whether attracting cycles might exist, or whether the substantial numerical evidence points to general convergence.

\subsection{Starting points; order of update of approximation vector}
As mentioned earlier, there is no theory yet about where to start the Ehrlich--Aberth-iteration from. For our experiments, we felt it was natural to use those points that provide good starting points for Newton's method: hence, we used a vector of $d$ points equidistributed on a large circle that surrounds all the roots. In contrast, Bini \cite{BiniAberth} discusses the use of Rouch\'e's theorem to find several circles that enclose different clusters of roots. 

The update from the ``old'' vector in $(z_1,\dots,z_d)\in\C^d$ of possible roots to the ``new'' vector $(z'_1,\dots,z'_d)$ can be done in (at least) two ways: for the computation of the new coordinates $z'_k$ one can either use only the old coordinates $(z_1,\dots,z_d)$ throughout (``Jacobi-style''), or, for each $z_k$, make use of the already computed new coordinates $z'_1,\dots,z'_{k-1}$ (``Gauss--Seidel-style''). The former approach is a natural iteration in $\C^d$, while the latter is not, but it is considered more efficient.

\subsection{The MPSolve Implementation} \label{Sec:MPSolveImpl}

There is a prominent software package that finds roots of univariate polynomials in practice, authored by Dario Bini, Giuseppe Fiorentino and Leonardo Robol and called \pack{MPSolve} (with the last published version \pack{3.1.8}). This software package has an impressive track record on root finding. Experiments were made for polynomials of degrees exceeding $25\,000$, some of them with coefficients larger than $10^{300}$; see \cite{BiniAberth,RobolMaster}. For all these experiments, usually no more than 20 Ehrlich--Aberth iteration cycles sufficed to find all roots (due to superlinear convergence, the required precision for root finding is not a decisive issue once all roots were found with sufficient precision so as to separate them).

In this implementation, the update of the new approximate root vectors is done in Gauss--Seidel style: every new coordinate $z'_k$ was used immediately in all subsequent computations.

\section{Newton's method and iterated refinement: virtues and problems}
\label{Sec:Newton}

Newton's method in its original form is a heuristic not an algorithm: it gives a rule to modify an initial guess, and it comes with a promise that initial guesses sufficiently close to a root will converge to this root. However, many starting points will not converge to a root, and if they do, it is not clear where to start in order to find all roots, not a subset.

For a polynomial $p$, the Newton method $N_p$ is a rational map that very naturally forms a dynamical system on the Riemann sphere $\Cbar$. The Riemann sphere decomposes into Julia set $J$ and Fatou set $F$, both of which are invariant under the dynamics. The Fatou set is open and dense and contains all the roots. The Julia set is compact and nowhere dense, but it is a non-empty, even uncountable compact set that contains finitely many periodic points of all periods $m\ge 2$. Any starting points in the Julia set will remain in the Julia set forever, so these points fail to converge to any root. It is known that the Julia set may have positive 2-dimensional Lebesgue measure, so a random starting point may be in the Julia set.

\begin{figure}[htbp]
\includegraphics[width=.6\textwidth]{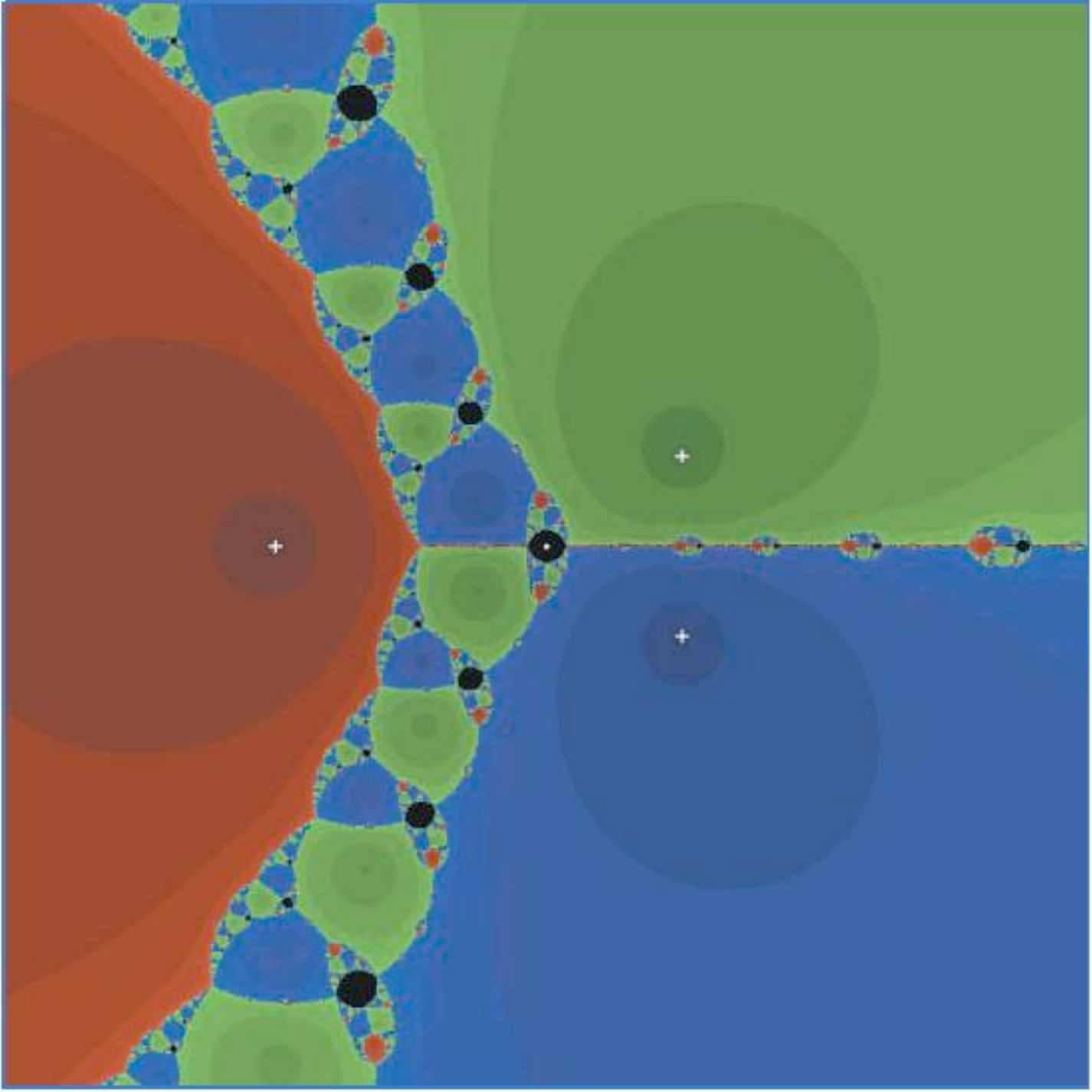}
\caption{The  Newton map for $p(z)=z^3-2z+2$ has a superattracting 2-cycle; the domain in black converges to this cycle. }
\label{Fig:NewtonCubicAttractingCycle}
\end{figure}

Worse yet, the Fatou set may contain open sets of points that converge to attracting cycles of any period $m\ge 2$ and hence not to any root. The simplest such example occurs for $p(z)=z^3-2z+2$ where there is a superattracting 2-cycle $0\mapsto -1\mapsto 0$ for the Newton map $N_p(z)=\frac{2z^3-2}{3z^2-2} $; see Figure~\ref{Fig:NewtonCubicAttractingCycle}, but it may happen for any degree and any period; a complete classification may be found in \cite{NewtonClassification}.

Therefore, individual orbits may fail to converge to a root (even though one may argue that most orbits will converge to a root; see \cite[Section~4]{NewtonInventiones}). If nothing is known about the positions of the roots (other than a normalization such as a circle enclosing all the roots), then one needs more than $d$ starting points in order to have a chance that all $d$ roots will be found. A sufficient such set of starting points has been constructed explicitly in \cite{NewtonInventiones}.

These starting points  have to be placed on a circle, or on a small number of circles, that surround all roots and are at a certain distance from the roots, at places where a good control of Newton's method can be assured. The disadvantage is that from there it takes a lot of iterations until the orbits reach the region where the roots are; typically, more than $d$ orbits each have to undergo more than $d$ iterations before interesting dynamics happens: so $O(d^2)$ iterations are ``wasted'' to move from the domain with control about Newton's method to the domain where the roots are. This prohibits any complexity of less than $O(d^2)$ iterations.

\begin{figure}[htbp]

{\includegraphics[width=0.7\textwidth]{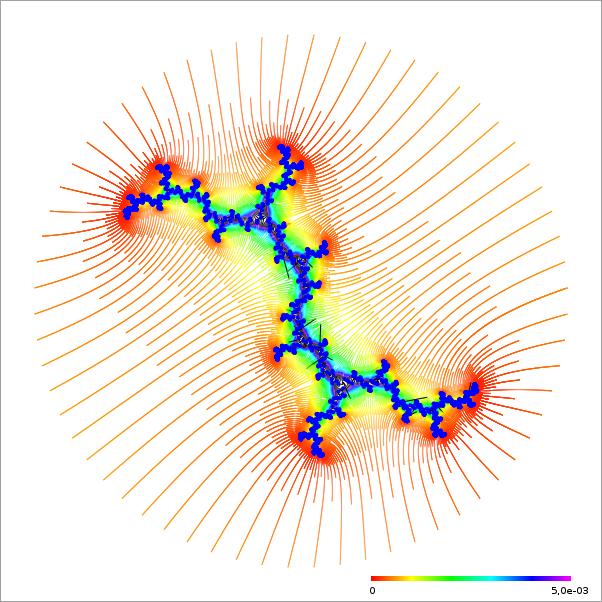}}
	\caption{Illustration of the ``iterated refinement'' idea: we start with a small number of points and refine when adjacent orbits no longer move in parallel. Note that most of the refinement happens once the orbits are close to the roots, which drastically reduces the necessary computations. Figure taken from \cite{NewtonRobin2}.}
\label{Fig:IteratedRefinement}
\end{figure}

In response to this problem, the ``iterated refinement'' Newton method has been developed in \cite{NewtonRobin1,NewtonRobin2}: the idea is that initially, where the orbits are still close to the initial circles, there is a lot of control on the Newton dynamics so all orbits do just about the same thing. The iterated refinement approach thus starts with a small initial number of starting points on the original circles, perhaps 64 points, and iterates them while their orbits are ``parallel'' in a sense to be made precise (for instance, in the sense of cross ratios between three adjacent orbits, together with $\infty$). One may also think about such orbits as of ``boring'', they don't exhibit any unexpected behavior. Once two or three orbits fail to remain parallel, these orbits are refined by creating additional orbits between them. We think of this occurrence as of indication that dynamics becomes more intricate and needs more orbits to be studied. This way, for certain families of polynomials of degrees up to $2^{30}>10^9$, all roots were found with log-linear complexity and log-linear computing time \cite{NewtonRobin2}. For an illustration of the iterated refinement in action one can have a look at Figure \ref{Fig:IteratedRefinement}.

The problem is that this is experimental science, and the original guarantee that all roots will be found by $d$ or more starting points is given up in favor of speed of computation. Among the heuristic parameters are the exact way to measure deviation of parallel orbits, as well as a quantitative refinement sensitivity with respect to this measure. Too insensitive refinements will lead to missed roots, while too high sensitivity will lead to slow computation. It is not clear how to fill in roots that are missed in the end. Deflation is usually not an option: it is in general very unstable, and it may destroy the simple structure of the polynomial (for instance, if the polynomial is given by iteration, then after deflation this structure is lost, and the evaluation becomes highly inefficient).

There is another possible optimization for Newton's method: for those cases where polynomials are given in terms of coefficients, one has numerous parallel Newton orbits, hence has to evaluate the same polynomial $p$ at many different points. We did these evaluations at all these points independently, not taking advantage of possible efficient evaluation of polynomials at many points in parallel. 

\hide{
===
experimental parameters, missing roots can not be recovered, roots often missed in parallel approach --- but also at exposed points???

parameters to be adjusted
* number of starting points (fix adjusted to degree?), sensitivity for refinement, criterion for refinement, circle radius for starting points, max number of refinement generations, max number of iterations, stopping criterion

postprocessing by theoretical guarantee (disk bounds) or by E--A
}


\section{Experimental results of both methods in comparison}\label{Sec:Experiments}

\subsection{Overview of polynomial families considered}

The numerical experiments were performed on various classes of polynomials of different degrees, defined in as different ways as possible, and so that the geometry of the root distribution was as different as possible, in order to explore different behavior of the two root finders. We explored the following families of polynomials:
\begin{enumerate}
\item
iterated quadratic polynomials (or periodic points of quadratic polynomials), Section~\ref{Sub:ExpIteratedQuadratic};
\item
polynomials that describe the centers of hyperbolic components of the Mandelbrot set, Section~\ref{Sub:ExpMandelbrotCenters};
\item
Chebyshev polynomials, Section~\ref{Sub:ExpChebyshev};
\item
Legendre polynomials, Section~\ref{Sub:ExpLegendre};
\item
polynomials with random roots on a circle, Section~\ref{Sub:ExpRandomOnCircle};
\item
polynomials with random roots on a disk, Section~\ref{Sub:ExpRandomInDisc}; 
\item
polynomials with roots on a finite square grid, Section~\ref{Sub:ExpGrid}.
\hide{rectangle = square?}
\hide{
\item
polynomials with random roots on a semi-circle;
\reminder{omitted?}
\item
polynomials with random roots in a neighborhood of a circle;
\reminder{omitted?}
\item
polynomials where the roots form several clusters (near-multiple roots);
\reminder{omitted?}
\item
polynomials where the roots are randomly distributed on a segment.
\reminder{Same as interval? Omitted?}
}
\end{enumerate}

The motivation for these various cases was quite different: iterated polynomials were originally chosen because they can be evaluated easily, and our focus was on root finding not polynomial evaluation. Chebyshev and Legendre polynomials were selected because the distribution of their roots is well known (along intervals). For polynomials where the coefficients are equidistributed \cite{ErdoesTuran,RootDistributionArnold}, the roots tend to accumulate at the unit circle with even distribution; this justifies some of the other choices.  

\subsection{Evaluation of polynomials}

A polynomial of degree $d$ is a holomorphic mapping from $\C$ to $\C$ of degree $d$; we would like to stress that one should not assume by default that it is given in terms of coefficients: indeed, all we require is that $p$ and $p'$ can be evaluated in some way (or possibly only $p'/p$). Indeed, much of the difficulties in finding roots of polynomials lies in evaluating the polynomials, and this problem is especially pronounced for high degree polynomials that are given in terms of coefficients. Here is one of our favorite prominent examples: periodic points of period $n$ for $p(z)=z^2+2$ are roots of the polynomial $P_n(z)=p^{\circ n}(z)-z$, where $p^{\circ n}$ denotes the $n$-th iterate of $p$. The constant coefficient in $P_n$ has size greater than $2^{2^n}$, which is utterly impossible to handle even for small $n$, while $P_n$ and its derivatives can easily be evaluated for $n$ as large as $25$ or more \cite{NewtonRobin2}, and all these $2^{n}$ roots have been found successfully and easily. 

In our experiments, some polynomials are given by recursive formulas that allow for evaluation in logarithmic time with respect to the degree. Others are evaluated in terms of coefficients, and yet others are, for the purposes of these experiments, evaluated in terms of the roots that are known ahead of time: in some experiments we had wanted to use a prescribed distribution of the locations of the roots, so these roots were computed first and the polynomials were evaluated in terms of these. This may be against the spirit of ``finding'' the roots, but is quite interesting for checking properties of root finding algorithms. Our main focus is not on polynomial evaluation, but we do want to point out that the evaluation is a serious and interesting issue, and often a serious bottleneck in the performance of root finders.

\subsection{Remarks on implementation and performance}

Our focus was on flexibility of experimentation, not real time optimization. We thus did most of our implementations in \pack{Java}.

We measured complexity in terms of number of arithmetic operations, and that is independent of the environment of implementation. More precisely, our implementation contains a counter which is updated each time an operation of real addition or multiplication is performed. 

\hide{In many cases, the practical limitations were in terms of evaluation of large polynomials, not computing time.} In order to work with high-degree polynomials we had to solve potential overflow problems. In particular, when we compute any fraction where numerator or denominator can be unexpectedly large we manually maintain the scientific representation of numbers by remembering mantissa and the order of magnitude.

In cases when a polynomial is given by its roots we computed the value of $p(z)/p'(z)$ with the formula $\left( \sum_{i=1}^d \frac{1}{z-\xi_i} \right)^{-1}$.

When performing the Ehrlich--Aberth iteration we update the components of the approximation vector in Gauss-Seidel style. However all components are updated until the iteration stops, even if some components approximate the root well enough. We admit that our implementation of the Ehrlich--Aberth iteration was quite down-to-earth. We believe, however, that the overall structure of the results is unaffected by such possible optimizations. A relevant method of improvement in the special case of polynomials with ``fast'' evaluation was recently pointed out to us by Dario Bini; this is discussed in Section~\ref{Sub:Improvements}.

\subsection{Experimental results}

We describe the results of our experiments on the performance comparison. Graphs picture the dependence of number of operations (real additions and multiplications) needed to find all the roots in terms of the degree of polynomials. As mentioned earlier, various polynomials were evaluated in linear time with respect to the degree (``slow evaluation'') and others in logarithmic time (``fast evaluation''), depending on their definition.


\subsubsection{Iterated quadratic polynomials}
\label{Sub:ExpIteratedQuadratic}

A class of polynomials that is particularly easy to evaluate are those that describe periodic points of period $n$ for quadratic polynomials $p_c(z)=z^2+c$, with $c\in\C$: these polynomials have the form $P_{n,c}(z)=p_c^{\circ n}(z)-z$, where $p_c^{\circ n}$ denotes the $n$-th iterate of $p_c$. These polynomials have degree $2^n$ but can be evaluated in complexity $O(n)$ hence in logarithmic complexity with respect to the degree. For certain choices of $c$, all periodic points of periods up to $30$ (i.e., degrees up to $2^{30}>10^9$) were computed successfully by Newton's method in \cite{NewtonRobin1,NewtonRobin2}.

\begin{figure}[htbp]
\framebox{
\includegraphics[width=.9\textwidth]{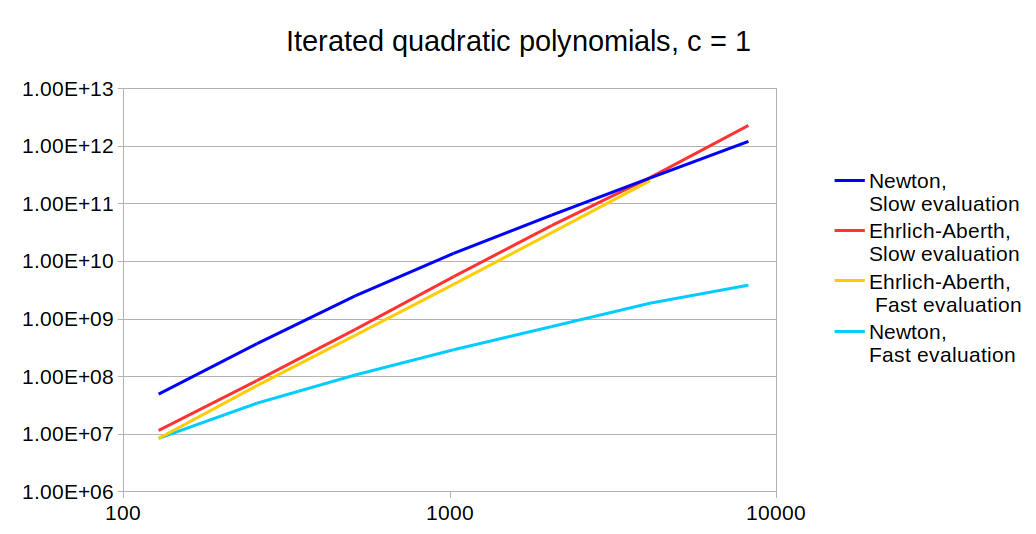}}
\framebox{
\includegraphics[width=.9\textwidth]{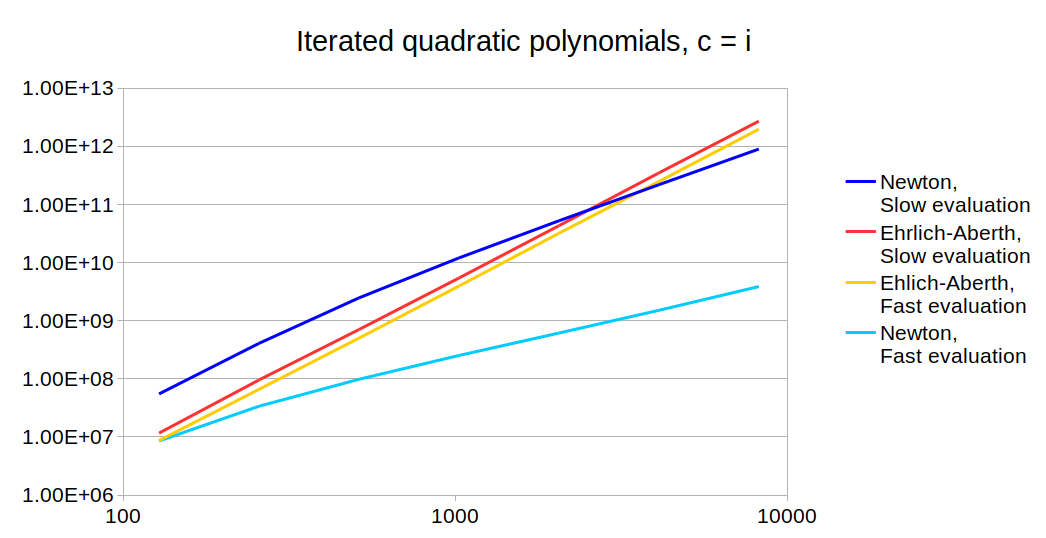}}
\caption{Complexity  (in terms of arithmetic operations) for finding periodic points of quadratic polynomials; log-log scale. Top: the  polynomial $z^2+1$; bottom: the polynomial $z^2+i$.
This graphs combine results for evaluating the polynomials in coefficient form (``slow evaluation'') and evaluating using iteration (``fast evaluation''). Notice that for the Ehrlich--Aberth-iteration does not accelerate significantly by taking advantage of fast evaluation.
}
\label{Fig:IteratedQuadratics}
\end{figure}

In order to gauge the dependence of complexity of the possibility to evaluate the polynomials efficiently, these experiments were run twice: once by evaluating $P_{n,c}$ in logarithmic time (``fast'') and once in linear time (``slow''). The results are shown in Figure~\ref{Fig:IteratedQuadratics}. It turns out that Newton's method is much faster, especially for high degrees, when taking advantage of the fast polynomial iteration, but not \reminder{necessarily} otherwise. Moreover, the  graphs show clearly that most of the computing time for Ehrlich--Aberth is not spent for the evaluation of the polynomials (but compare the discussion in Section~\ref{Sub:Improvements}).

It is a well known and simple fact from complex dynamics that periodic points of $p_c$ are distributed along the Julia sets (even equidistributed with respect to harmonic measure). For the later discussion it may be worth mentioning that the Julia sets are never convex, except when $c=0$ (the Julia set is a circle) or $c=-2$ (the Julia set is an interval). Our choice of polynomials are $z^2+1$ (the  Julia set is a Cantor set) and $z^2+i$ (the Julia set is a dendrite).

\hide{
\reminder{more precisely: what polynomials, what was done? In the report it says, evaluating in terms of ther roots found earlier; that is unusual, why not using coefficients?}  \reminder{Show picture of root location! Once coefficients known.}
}



\subsubsection{Centers of hyperbolic components of the Mandelbrot set}
\label{Sub:ExpMandelbrotCenters}

\begin{figure}[htbp]
\framebox{
\includegraphics[width=.8\textwidth]{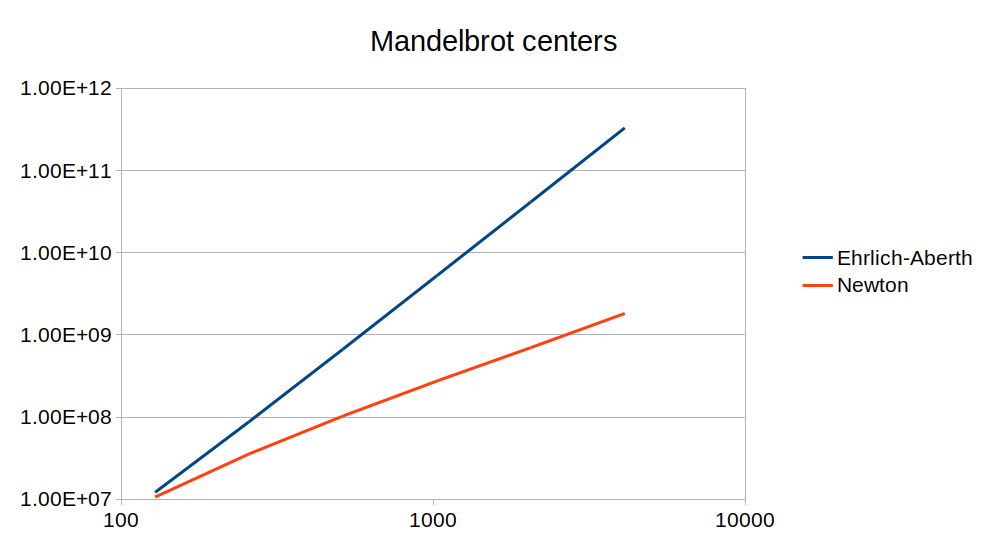}}
\caption{Complexity for finding centers of hyperbolic components of the Mandelbrot set; log-log scale. 
}
\label{Fig:Mandelbrot}
\end{figure}

These polynomials in the variable $c$ are defined by the iteration $q_0=0$ and $q_{n+1}=q_n^2+c$, i.e.\ $q_1=c$, $q_2=c^2+c$, $q_3=(c^2+c)^2+c$, etc., so $q_n$ has degree $2^{n-1}$. Equivalently, for the iteration $p_c\colon z\mapsto z^2+c$, we have $q(c)=p_c^{\circ n}(z)$ with $z=0$: the roots are exactly the parameters $c$ for which the  critical point $z=0$ of $p_c$ is periodic with period $n$. Such parameters are well known to be the centers  of hyperbolic components of period $n$ for the Mandelbrot set (see for instance \cite{Orsay,DS_Dynamics}). These polynomials can be evaluated by recursion as efficiently as the iterated quadratic polynomials from Section~\ref{Sub:ExpIteratedQuadratic}. For large $n$, these points accumulate at the boundary of the Mandelbrot set, so they form a very non-convex set. The results of the experiments are displayed in Figure~\ref{Fig:Mandelbrot}.

It may be interesting to note that these polynomials have been used as sample polynomials for evaluating the efficiency the Ehrlich--Abert-method  in its implementation in  \texttt{MPSolve} (see \cite{MPSolveExp}), and independently also for our experiments on Newton's method (see \cite{NewtonRobin1,NewtonRobin2}). One has to keep in mind that formally speaking, \cite{MPSolveExp} works with polynomials $q_n(c)/c$ of degree $2^{n-1}-1$ but this difference is, of course, insignificant for root finding.

We recently learned that a new implementation of Ehrlich--Aberth's method that focuses especially on polynomials given by iteration  has yielded particularly efficient evaluation for high degrees; see the discussion in Section~\ref{Sub:Improvements}.

\subsubsection{Chebyshev polynomials}
\label{Sub:ExpChebyshev}

These polynomials have the special property that they can be evaluated in logarithmic time with respect to the degree, based on the recursion formula $T_{2d}(x)=2T_d(x)^2-1$. This evaluation works, of course, only for degrees $d = 2^n$, but we have to mention that Chebyshev polynomials can be defined for arbitrary degree. Two sets of experiments were done for these polynomials as well, with evaluation in logarithmic and in linear time (``fast'' and ``slow''). 

For these polynomials, Newton turned out to be significantly fast than Ehrlich--Aberth in both cases --- obviously more so for the recursive evaluation.

\begin{figure}[htbp]
\framebox{
\includegraphics[width=.8\textwidth]{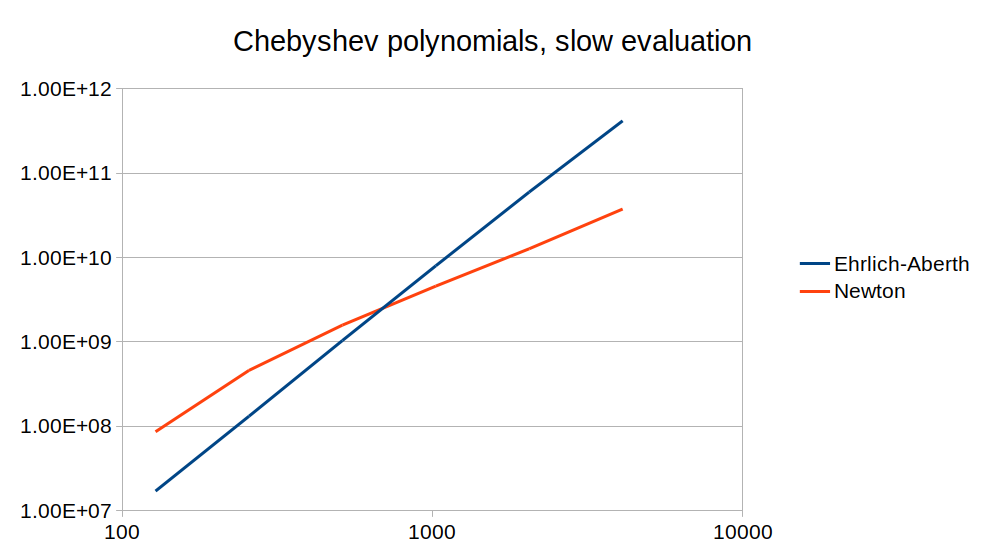}}
\framebox{
\includegraphics[width=.8\textwidth]{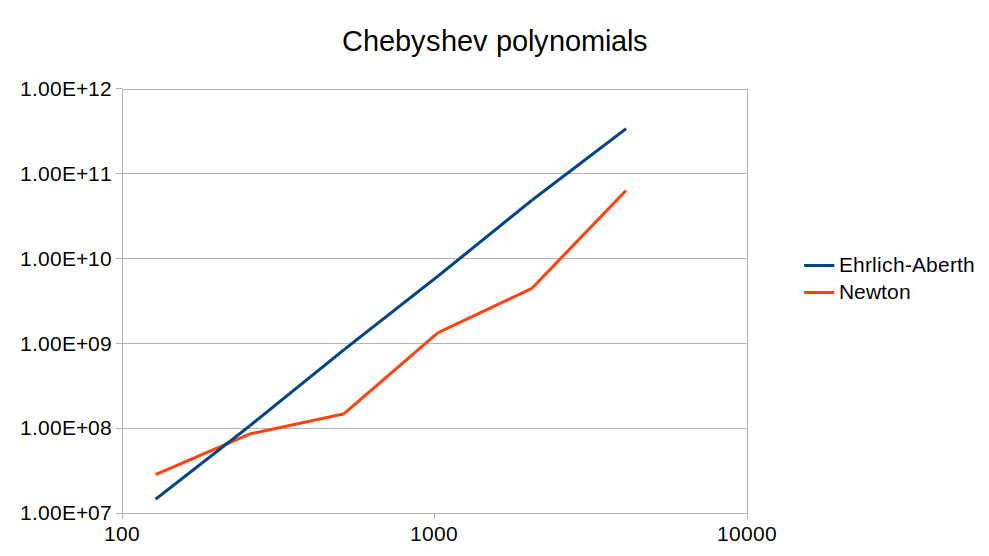}}
\caption{Complexity for finding periodic points of Chebyshev polynomials; log-log scale. Top: evaluation of polynomials in coefficient form (``slow evaluation''); bottom: evaluation using iteration.
}
\label{Fig:Chebyshev}
\end{figure}


\subsubsection{Legendre polynomials}
\label{Sub:ExpLegendre}

Here the evaluation of the polynomials was always in linear time. The result of the observed complexity seems comparable to the ``slow'' experiment of Chebyshev polynomials (Figure~\ref{Fig:Chebyshev}); which is consistent with the observation that the speed of evaluation is comparable between both polynomial families, and so is the location of the roots (on the unit interval).

\begin{figure}[htbp]
\framebox{
\includegraphics[width=.8\textwidth]{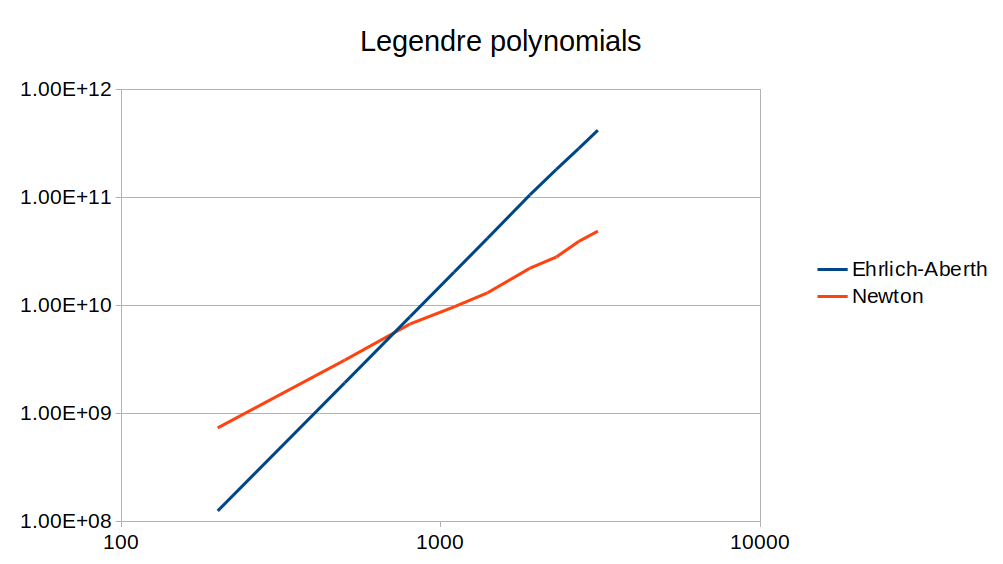}}
\caption{Complexity for finding roots of Legendre polynomials; log-log scale. Evaluated ``slowly''.
}
\label{Fig:Legendre}
\end{figure}


\subsubsection{Polynomials with random roots on a circle}
\label{Sub:ExpRandomOnCircle}

Here the structure of the experiment was that initially a given number of roots was distributed randomly on a circle, and then the polynomials were defined in terms of these roots (evaluation in linear complexity). The outcome is that Newton's method seem to be faster for large degrees, but this effect is not very pronounced.
\begin{figure}[htbp]
\framebox{
\includegraphics[width=.8\textwidth]{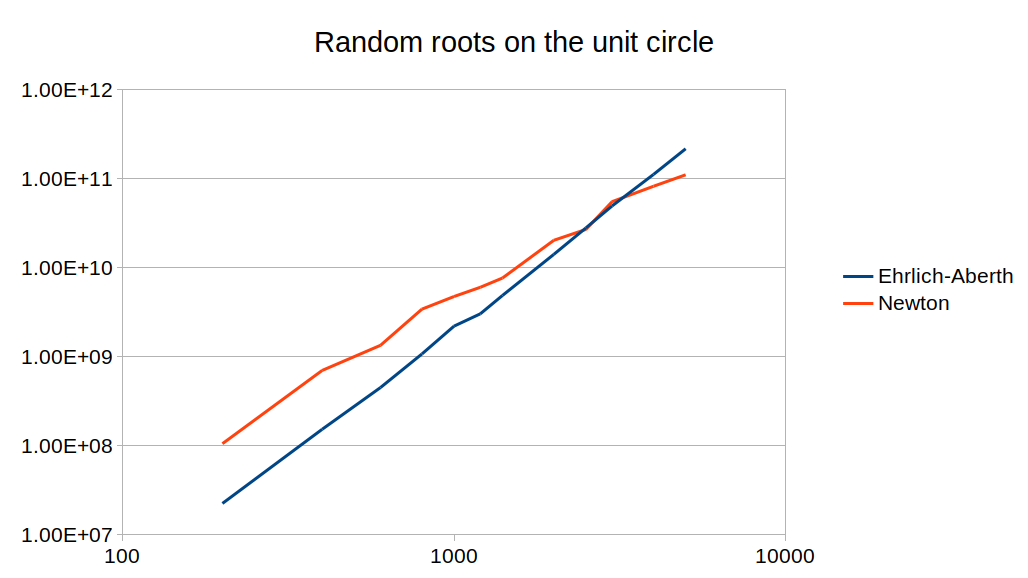}}
\caption{Complexity for finding roots distributed randomly on the unit circle; log-log scale.
}
\label{Fig:RandomOnCircle}
\end{figure}

\hide{
\reminder{In the experiments with roots on a unit circle the graph has a weird jump in the end. Omit and postpone for next round when this jump is clarified?}
}


\subsubsection{Polynomials with random roots in a disk} \label{Sub:ExpRandomInDisc}

For these polynomials, the roots were placed randomly in the unit disk, and the polynomials were again defined in terms of these roots; as a result, their evaluation is linear in the degree (``slow''). In this case, the Ehrlich--Aberth method is consistently faster.

We should mention that in our implementation, the roots were distributed independently over the unit disk, but not with respect to planar Lebesgue measure, but with respect to the linear product measure in polar coordinates: i.e. for the points $z=re^{i\phi}$, the distribution used planar Lebesgue measure for $(r,\phi)\in[0,1]\times[0,2\pi]$.

\begin{figure}[htbp]
\framebox{
\includegraphics[width=.8\textwidth]{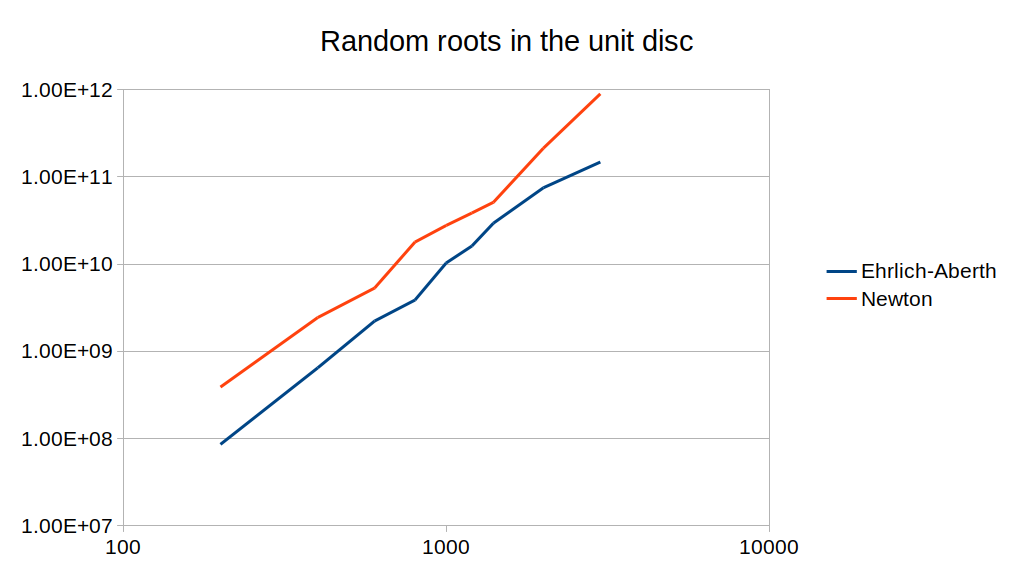}}
\caption{Complexity for finding roots distributed randomly in the unit disk; log-log scale.
}
\label{Fig:RandomInDisc}
\end{figure}


\subsubsection{Polynomials with roots on a square grid} 
\label{Sub:ExpGrid}

In this set of experiments, the roots were placed on a grid $\{1,2,\dots,n\}+ i\{1,2,\dots,n\}$ of points in $\C$, for various values of $n\in\{4,40\}$, and once again the polynomials were defined in terms of their roots, with evaluation in linear time. Consistently, the Ehrlich--Aberth method was faster than Newton's method, but not as significantly as for random roots on the disk.

\begin{figure}[htbp]
\framebox{
\includegraphics[width=.8\textwidth]{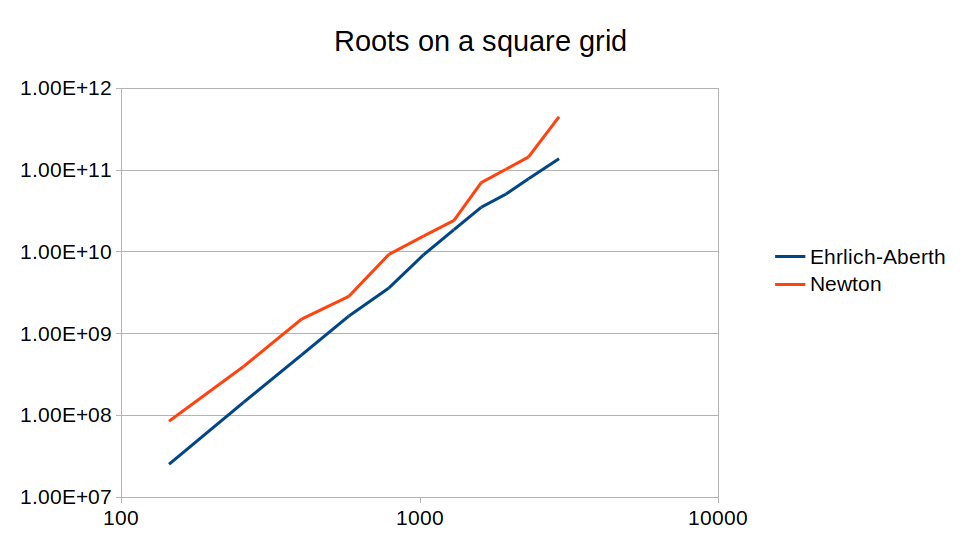}}
\caption{Complexity for finding roots on a square grid; log-log scale.
}
\label{Fig:Grid}
\end{figure}

\hide{
\sergey{maybe also add random on a line in a disc and point out the similarity to slow Chebyshev and Legendre: geometry of roots is more important than particular properties of a polynomial}
}


\hide{
\reminder{NOTE: discuss complexity: quadratic for every EA iteration step, possibly with small number of iterations; Newton can be parallelized / parallel evaluation, we have not that; good overall complexity for Marvin experiments.}
}

\section{Conclusion}
\label{Sec:Conclusion}

\subsection{Analysis of the performance}

The outcome of our experiments is that, among our two implementations, none of them is universally better: for some cases, one of the two methods outperforms the other, and for others the result is the opposite. What follows is a heuristic explanation for the observed behavior.
\begin{enumerate}
\item \emph{Evaluation of polynomials.}
Some of the experiments were motivated by the choice of polynomials that can be evaluated in logarithmic time with respect to the degree. This leads to efficient evaluation by Newton's method. Newton's method generally outperforms Ehrlich--Aberth for those polynomials: this includes iterated quadratic polynomials (Section~\ref{Sub:ExpIteratedQuadratic}) and Chebyshev polynomials (Section~\ref{Sub:ExpChebyshev}), as well as the centers of hyperbolic components that are among the ``demo polynomials'' for the MPSolve implementation. 

This advantage plays an especially significant role when the degrees of the polynomials are large, and for these has the decisive effect on the overall efficiency outcome.

\item \emph{Location of roots.}
There are several additional cases when Newton outperforms Ehrlich--Aberth, even though the evaluation of the polynomials was only linear with respect to the degree: among our experiments these were Chebyshev polynomials (Section~\ref{Sub:ExpChebyshev}), Legendre polynomials (Section~\ref{Sub:ExpLegendre}), and those polynomials where the roots are randomly distributed on the unit circle (Section~\ref{Sub:ExpRandomOnCircle}). These were exactly those experiments where the roots were ``exposed'' from infinity, that is all roots were on the boundary of the convex hull formed by the roots

On the other hand, Ehrlich--Aberth was faster, sometimes very much so, when roots were ``hidden'' in the interior of the convex hull: this includes the case of iterated quadratics with linear (``slow'') evaluation (Section~\ref{Sub:ExpIteratedQuadratic}) where the roots are distributed along polynomial sets that are almost never convex, as well as the cases of random roots on the disk (Section~\ref{Sub:ExpRandomInDisc}) and when the roots were distributed on a grid (Section~\ref{Sub:ExpGrid}).

\hide{
\reminder{Have you counted how many E--A iterations were usually required? I think I remember Bini saying that overall about six iterations were sufficient. Perhaps that calls for systematic experiments too?}
}

\hide{	
	The difference between these polynomials is in the root configuration. Polynomials from items \ref{Sub:ExpChebyshev}, \ref{Sub:ExpLegendre} and \ref{Sub:ExpRandomOnCircle} has more ,,exposed'' configuration of roots \sergey{better seen from infinity/ starting points of iteration; on a boundary of convex hull of roots.} Polynomials from items \ref{Sub:ExpRandomInDisc} and \ref{Sub:ExpGrid} has most of the roots ,,hidden'' by the outer roots.

		The most convincing explanation to us is that the Newton's iterated refinement method starts from the neighborhood of infinity and moves towards the roots. Other roots on its way create an obstacle where iterations converge not that fast. The geometry of E--A method is different and iterations don't approach roots only ,,from outside''.
}
\end{enumerate}

\subsection{Complexity}

We are not aware of any complete analysis of the complexity of both methods, but there are some partial results that shed some light on the overall efficiency.

For the Ehrlich--Aberth-method, the complexity of every iteration step seems to be of the order $O(d^2)$, irrespective of how efficiently the polynomials $p'/p$ can be evaluated. In the various experiments performed by Bini and Robol, no more than 20 iteration steps needed to be executed overall in order to find all roots (where the required precision for finding the roots is largely irrelevant, at least in case of simple roots, because of the superlinear convergence once all roots are separated). This suggests an overall complexity of $O(d^2$) for finding all roots. It is not sure that this method is generally convergent, even though failure to converge has never been observed experimentally except in cases with clear symmetry. If the method was not generally convergent, then this may have a relevant impact on the complexity.

In the case of Newton's method, the results are less clear, especially since the method is not generally convergent, and orbits that fail to converge can frequently be observed. One needs at least $d$ orbits to iterate, and can be certain to find all roots when starting at $O(d\log^2d)$ initial points  \cite{NewtonInventiones}. Each of these orbits starts away from the disk containing the roots, and in order to have good control on the starting points need to be iterated at least $d$ times each in order to land in the disk containing all roots; this gives a lower bound on $O(d^2\log^2d)$ Newton iterations. It is shown in \cite{NewtonTodor} that one can expect, under reasonable assumption of equidistribution of the roots, that the same complexity should also suffice to find all roots. Of course, this complexity needs to be multiplied with the complexity to evaluate each Newton step, which may be $O(d)$ or $O(\log d)$ depending on how the polynomial is specified. However, parallel evaluation of all the Newton orbits can yield an overall logarithmic complexity for the Newton evaluation, so that we expect again arithmetic complexity of $O(d^2\log^kd)$ for some $k$ (based on a machine model that has infinite precision in each operation); however, there may be issues with stability.

In \cite{NewtonRobin2}, the iterated refinement Newton method was introduced; this was also the procedure implemented here as described earlier. The number of orbits that are iterated in parallel goes down from $d\log^2d$ to a much smaller number. The exact complexity depends on the sensitivity of refinement and is difficult to estimate theoretically (too sensitive refinement leads to $d$ orbits that are iterated in parallel, hence no gain in complexity; too insensitive refinement leads to missed roots that are hard to recover.)

\subsection{A new conjecture on starting points for Newton's method}

Neither of the two root finding methods discussed here has a complete theory about its global dynamics. Very little is known about the global properties of the Ehrlich--Aberth-method (the recent results in \cite{BernhardMichaelWeierstrass} show that the related Weierstrass method is not generally convergent, but it is not clear what the conclusion about Ehrlich--Aberth should be). More is known about Newton's method, beyond the fact that  it is not generally convergent and may have attracting cycles of any period. The first quantitative result was published in \cite{NewtonInventiones}, as mentioned earlier; it is usually presented for polynomials with all roots in the unit disk $\disk$ (this is not a restriction after an appropriate coordinate change).

\begin{theorem}[Starting points for Newton's method, \cite{NewtonInventiones}] 
For every degree $d$, there are $c_1\log d$ circles with $c_2d\log d$ points on each so that for every  polynomial $p$ of degree $d$  with all roots in  $\disk$, for every root $\alpha$ of $p$ at least one of these $c_1c_2d\log^2d$ points converges to $\alpha$ under iteration of the Newton method.
\end{theorem}

These circles and the points on them, as well as the constants $c_1$ and $c_2$ are explicitly given. When allowing to take the circles far away from $\disk$ (thus allowing many iterations until the points on the circle reach the disk containing  the roots), one may get $c_1c_2<1.1$.

This result was improved in \cite{NewtonProbabilistic}, using a probabilistic set of starting points, to a total number of starting points of size $c_3d(\log\log d)^2$ that finds all roots with arbitrary desired probability (where of course the constant $c_3$ depends on this probability). 

In the experiments described in this paper, we discovered that Newton's method finds roots more easily when they are on the boundary of the convex hull of all the roots. We believe that the complexity to find all roots with this property is much better, and suggest the following.

\begin{conjecture}
For every radius $r>1$ there is a universal constant $c>1$ with the following property. If $p$ is a polynomial of degree $d$, rescaled so that all its roots are in the unit disk, then $\lceil c d\rceil$ equidistributed points on the circle $\partial D_{r}(0)$ have the property that for each root on the boundary of the convex hull of all roots, at least one of the given points will converge to this root under  iteration of Newton's method. 
\end{conjecture}

In other words, there is a universal set of $O(d)$ points that finds all the $d$ roots of $p$ that are on the boundary of the convex hull. We believe that $c\approx 2.6$ for big $r$.  

This conjecture goes back to the Bachelor thesis \cite{Sergey_Bachelor}. In order to express the underlying idea, we need to review some background from \cite{NewtonInventiones}. For every root $\alpha$, the \emph{immediate basin} of the root is the connected component containing $\alpha$ of the set of points that converge to $\alpha$ under the Newton iteration. A \emph{channel} is a homotopy class of curves in $U_\alpha$ with endpoints fixed that connect the root $\alpha$ to $\infty$. Every root has some number $k_\alpha$ of channels with $1\le k_\alpha\le d-1$. Each of these channels has a certain ``thickness'', measured in terms of a conformal invariant. The thicker a channel, the fewer starting points are needed in order to be sure that one of them is in the immediate basin. The more channels a root has, the thinner they may be, and the more starting points are required.

The idea behind the conjecture is that every root on the boundary of the convex hull has ``main channel'' with the property that the root can be connected to $\infty$ along this main channel so that the convex hull of the roots is avoided. While this is a theorem, the conjecture is based on the assumption that a main channel has at least the same thickness as a root with only two channels. The conjecture would then follow as in \cite{NewtonInventiones} in the case when all roots are real.

Of course, with the conjectured improvement of the number of starting points, the complexity of finding all the roots on the boundary of the convex hull, as described in \cite{NewtonEfficient,NewtonTodor} would also improve.

\subsection{Postprocessing; mixing both approaches}

One notorious issue for the Iterated Refinement Newton Method is that it is optimized beyond the regime for which convergence can be guaranteed by theoretical bounds; it depends on some experimental parameters such as the refinement threshold. In case of too optimistic refinement, hence too insensitive refinement threshold, a certain number of roots may be missed in the process. One possible conclusion is to run a ``postprocessing step'' of Ehrlich--Aberth on the approximate roots found by Newton. If most roots were found by Newton initially, say only $d'$ of the $d$ roots are still missing, then the complexity of each Ehrlich--Aberth-postprocessing step requires $O(d'd)$ operations, so this implies a nontrivial lower bound on the overall complexity. However, since the total number of Ehrlich--Aberth-steps is usually quite small (as mentioned in Section \ref{Sec:MPSolveImpl} and also confirmed by our experiments), this may become a realistic option only when $d'=o(d)$. 

This post-processing was implemented in various of our experiments. 
The general observation was that, compared to necessary post-pro\-ces\-sing by Ehrlich--Aberth, the performance of Newton's method was generally better when more sensitive refinement parameters were chosen that avoided any post-processing steps, even though this led to earlier refinement overall.

Another rather radical approach is to run both methods in parallel, sharing execution time equally between both processes, and stopping once one of the two methods terminates successfully. This seems wasteful as roughly half the operations are not used at all eventually. However, the overall complexity is then the minimum of both methods, and if the complexities of both of them differ significantly in ways that cannot be predicted efficiently ahead of time, this may lead to overall improved complexity. (After all, for each of the two methods our focus was on complexity not on improvement by bounded factors, so factors or two were often ignored.)

\subsection{Possible improvements}
\label{Sub:Improvements}

We take it as an encouraging sign that since the times of our original experiments, there have been inspiring discussions and suggestions for improvements on both methods.

Newton's method has its most obvious bottleneck in the case of polynomials with ``slow'' evaluation, i.e. when the computation of $p'/p$ has complexity $O(d)$. Since this is a  parallel computation, it can in principle be accelerated by computing all orbits in parallel. We were aware of fast algorithms that were, however, too unstable to be used in practice. We are grateful to Dan Spielman for having pointed out to us that meanwhile there are fast stable methods (albeit apparently not easy to implement). These should have the potential to speed up Newton's method even for general polynomials to a speed comparable to those with ``fast'' evaluation.

On the other hand, the Ehrlich--Aberth iteration step, especially in the form given in \eqref{Eq:EhrlichAberthIteration}, has two major ingredients: an evaluation of the logarithmic derivative $p'/p$, and the correction term 
$
\sum_{i\neq k}\frac{1}{z_k-z_i}
$.
The logarithmic derivative benefit from the same advantages of fast computation as for Newton's method (especially in the cases of particularly  suitable polynomials like the ones considered in some of our experiments, but also for general polynomials as described above). However, since the correction term requires $O(d^2)$ operations, its exact computation can far outweigh all computations depending on the polynomials (which is why we had not optimized  polynomial evaluations for the Ehrlich--Aberth-method). Soon after the first version of this paper was submitted, we learned from Dario Bini and Victor Pan about the Fast Multipole Method (FMM; see \cite{FMM}) that makes a very fast approximate computation of the correction term possible. Since all described methods only compute (hopefully ever-improving) approximations, it is not a structural disadvantage to settle for ``approximate approximations'' --- even more so as there is still no global theory for the Ehrlich--Aberth-method even in its precise form. 

Dario Bini informed us that he has used this improved Ehrlich--Aberth-method very recently for certain polynomials that were defined by fast recursion (possibly inspired by the super-fast performance of Newton's method for such polynomials) with almost linear CPU time up to very high degrees. 

It would be an interesting project to perform a fresh set of experiments where both methods are optimized in the ways described, and see whether the relative advantages or comparisons in terms of recursive form of the  polynomials, and shape of the roots, persist. This would, however, be an entirely new project, inspired  not the least by the results of the experiments described  here.

\subsection{Conclusion}

This project of comparing the efficiency between the Newton and Ehrlich--Aberth root finders is experimental, not an exact science. It is supposed to inspire or possibly provoke further investigations into the question which root finders perform more efficiently under which circumstances, and how this depends on the properties of the polynomials --- and possibly of the specific clever implementations of the root finders; and  we are quite pleased to learn that some of this has already happened since the first version of this paper was submitted. Surprisingly, even at this time the fundamental question of finding all roots of univariate complex polynomials seems far from resolved and invites for more  progress to be made!


\end{document}